\documentclass[]{article}

\usepackage [T2A] {fontenc}
\usepackage [utf8] {inputenc}
\usepackage[russian]{babel}

\usepackage{mathtools}
\mathtoolsset{showonlyrefs}

\usepackage{indentfirst} 
\usepackage[14pt]{extsizes}
\usepackage{amssymb}
\usepackage{amsmath}
\usepackage{amsthm}
\usepackage{pifont}
\usepackage{subfig}
\usepackage{graphicx}
\usepackage[colorinlistoftodos]{todonotes}
\usepackage{makecell}
\usepackage{multirow}
\usepackage{booktabs}
\usepackage[colorlinks=true, allcolors=blue]{hyperref}
\usepackage{algorithm}
\usepackage{algorithmic}

\def\R{\mathbb{R}}
\newcommand{\E}{{\mathbb E}}

\def\R{\mathbb R}
\def\E{\mathbb E}
\def\EE{\mathbb E}

\def\e{\varepsilon}

\def\ag#1{{\color{black}#1}}
\def\agv#1{{\color{red}#1}}
\def\gav#1{{\color{black}#1}}

\usepackage{geometry}
\def \EE {\mathbb E}

\def\e{\varepsilon}

\def \R {\mathbb R}

\begin{document}
\renewcommand{\abstractname}{\vspace{-\baselineskip}}

$$
\\\\
$$


\begin{center}
\Large{\textbf{Распределенные и параллельные алгоритмы решения задач анализа данных}}\footnote{Научно-популярный отчет, содержащий наработки по гранту РФФИ  19-31-51001 (Научное наставничество); Этап 1 и по работе, выполненной при поддержке Министерства науки и высшего образования Российской Федерации (госзадание) №075-00337-20-03, номер проекта 0714-2020-0005; Этап 1.}
\end{center}

\begin{center}
\textbf{Д.\,М.~Двинских, А.\,В.~Гасников,  \\ А.\,В.~Рогозин, А.\,Н.~Безносиков и др.}

\textit{
Московский физико-технический институт (национальный исследовательский университет), Долгопрудный, Мос. обл., Россия \\
}


\end{center}

\begin{abstract}
\noindent Одной из главных целей проекта была разработка оптимальных численных методов решения задач (вида суммы) децентрализованной распределенной (сильно) выпуклой оптимизации  с прямым и двойственным детерминированным и стохастическим оракулом (выдающим (стохастический) градиент прямой (двойственной) функции, находящейся в узле). Для распределенных алгоритмов существует два критерия, по которым можно измерять качество алгоритма: число обращений к оракулу в узле и число коммуникаций. Среди прочего, исследования, проведенные на первом этапе, показали, что 1) с точностью до логарифмических (по желаемой точности) множителей можно строить алгоритмы оптимальные по обоим критериям одновременно \cite{uribe2020dual,dvinskikh2021decentralized,rogozin2020towards,rogozin2021accelerated}; 2) При использовании прямого оракула в оценки сложности (вызовов оракула, число коммуникаций), на самом деле, должны входить средние (по узлам) константы гладкости и сильной выпуклости функций, находящихся в узлах, а не худшие, как казалось изначально. Более того, этот результат переносится и на коммуникационные графы, меняющиеся со временем \cite{rogozin2020towards,rogozin2021accelerated}; 3) можно строить теорию оптимальных безградиентных  численных методов решения задач децентрализованной распределенной (сильно) выпуклой оптимизации (вида суммы) с прямым детерминированным и стохастическим оракулом. Интересной это теория получается в негладком случае \cite{beznosikov2020derivative}.

\textbf{Ключевые слова}: выпуклая распределенная оптимизация, стохастический оракул, прямой оракул, двойственный оракул, оценка числа коммуникаций 
\end{abstract}

\section{Постановка задачи и \ag{основные полученные результаты}} \label{section_1}
\ag{В изложении материала данного раздела мы следуем в основном работе, опубликованной в журнале Ill posed Inverse problems \cite{dvinskikh2021decentralized}, см. также обзор \cite{gorbunov2020recent}.} 

В проекте рассматриваются преимущественно задачи выпуклой (для простоты считаем, что $f(x,\xi)$ -- выпуклая функция по $x$ для всех $\xi$) стохастической оптимизации следующего вида:
\begin{equation}
\label{SP}
\min_{x\in Q \subseteq \mathbb{R}^n} f(x) := \E[f(x,\xi)].    
\end{equation}

Основным источником таких задач в современных приложениях является анализ данных: математическая статистика \cite{spokoiny2012parametric} (получение оценки максимального правдоподобия -- истинное значение оцениваемого параметра является решением задачи стохастической оптимизации) и машинное обучение \cite{shalev2014understanding} (минимизация риска). Во втором случае закон распределения $\xi$ не известен. Можно только получать реализации $\{\xi^k\}_{k=1}^m$.   Параметр $m$ иногда называют объемом выборки. Требуется решить задачу \eqref{SP} с точностью $\e$ по функции, используя наименьший объем выборки $m$. Режим получения реализаций $\xi$ онлайн (в ходе работы алгоритма обучения) так и называют \textit{онлайн} подход. Противоположный подход (\textit{офлайн}) предполагает, что данные (выборка) уже где-то хранятся. Подробнее о том, что написано в данном абзаце можно посмотреть в \ag{лекции Александра Гасникова в рамках курса <<Оптимизация в машинном обучении>>, который читается участниками гранта (Э. Горбуновым, А. Рогозиным и др.) магистрам МФТИ в осеннем семестре 2021/2022 гг.:}
\begin{center}
    \small{\url{https://www.youtube.com/playlist?list=PLIvQImOQgbGZH-HlEsVYddBF6EU-qrDOv}.}
\end{center}

Далее в тексте мы будем многократно встречаться с распределенными (централизованными и децентрализованными) алгоритмами. Во всех этих случаях, следуя \cite{arjevani2015communication}, будем предполагать, что одно сообщение (акт коммуникации между узлами) имеет ограниченный объем $O(n)$. Для простоты мы будем измерять объем просто в числе передаваемых чисел выбранного типа. В действительности, тут можно говорить и о битах. С точностью до логарифмических поправок результаты останутся прежними \cite{alistarh2020improved}.

Все обсуждаемые далее способы параллелизации можно также понимать как централизованную архитектуру (только с общей памятью), в которой центральный узел (процессор) общается со всеми другими узлами (ядрами) также сообщениями размера (объема) не более $O(n)$. Во всех этих подходах в промежутках между актами коммуникации среди узлов осуществляются однотипные действия, которые в случае вычислительной эквивалентности узлов, приводят к одним и тем же временным потерям. Таким образом вопросы асинхронности вычислений здесь не рассматриваются. Под возможностью (полной) параллелизации работы метода на $p$ узлах понимается такой способ организации вычислений и обмена информацией, который приводит к сокращению времени работы метода (в приближении пренебрежения временем на коммуникации -- параллельный сбор и распространение данных / результатов вычислений) в $\sim p$ раз с сохранением общей (оракульной) сложности, т.е. с сохранением общего (суммарного) объема вычислений, который необходимо произвести.

\subsection{Негладкий случай} \label{section_1.1}
В классе онлайн алгоритмов при минимальных предположениях об $f$: 
\begin{equation}
\label{M}
\E[\|\nabla f(x,\xi)\|_2^2]\le M^2
\end{equation}
на итерациях алгоритма (т.е. при $(x,\xi) = (x^k,\xi^k)$), наилучшими (с точки зрения оракульных нижних оценок \cite{nemirovsky1983problem}) способом решения задачи \eqref{SP} в указанном выше смысле будет использование стохастического градиентного спуска SGD \cite{juditsky2012first-order,koloskova2020unified} (для простоты записи приводим SGD для $Q=\R^n$): 
$$x^{k+1} = x^k - h_k\nabla f(x^k,\xi^k)$$
 или его вариаций \cite{li2020root}, в зависимости от того, что выбирается в качестве критерия оценки сходимости. SGD дает следующую оценку числа итераций (объема выборки):
\begin{equation}\label{SNS}
m = m(\e) = \min\left\{O\left(\frac{M^2R^2}{\e^2}\right),~ O\left(\frac{M^2}{\mu\e}\right)\right\},
\end{equation}
необходимых для достижения точности $\varepsilon$ по функции в среднем.
При этом на каждой итерации требуется один раз посчитать $\nabla f(x,\xi)$. В приведенной оценке $\mu \ge 0$ -- константа сильной выпуклости $f$ в 2-норме,  а $\|x^0 - x^*\|_2 \le R$ -- расстояние от точки старта алгоритма  $x^0$ до решения $x^*$.

При офлайн подходе   задача \eqref{SP} заменяется на задачу \cite{shalev2009stochastic}:
\begin{equation}\label{empiric}
\min_{x\in Q \subseteq \mathbb{R}^n} \frac{1}{m}\sum_{k=1}^m f(x,\xi^k) + \frac{\e}{2R^2}\|x - x^0\|_2^2,
\end{equation}
с параметром $m$, который определяется с точностью до логарифмических множителей формулой \eqref{SNS}. В свою очередь, задачу \eqref{empiric} необходимо решать с точностью по функции $\e' = O(\max\{\mu,\e/R^2\}\e^2/M^2)$, чтобы $\e'$-решение \eqref{empiric} было $\e$-решением исходной задачи \eqref{SP}.

Отсюда видно, что офлайн подход соответствует онлайн подходу в смысле объема выборки, но в вычислительном плане выглядит менее привлекательным, если не учитывать возможность параллельных и распределенных вычислений. Собственно, так и было принято считать \cite{shapiro2014lectures}. 

\ag{Однако в начале 2020 года участником проекта Дариной Двинских было показано, что даже без возможности параллельных и распределенных вычислений, офлайн подход может работать быстрее онлайн \cite{dvinskikh2020sa}. Статья подана в журнал Optimization Methods and Software. В качестве примера задачи, на которой это удалось получить была выбрана задача вычисления барицентра Васерштейна вероятностных мер. В основе идеи Д. Двинских лежит наблюдение, что в данном конкретном случае оказывается выгоднее (по времени работы алгоритма) построить двойственную задачу к \eqref{empiric} и решать ее, чем осуществлять честную процедуру SGD-типа с вычислением стохастического градиента. Стоит отметить, что решать двойственную задачу можно и распределенными алгоритмам, тогда выгода будет только усиливаться \cite{gorbunov2019optimal}. Стоит также отметить совсем свежий результат \cite{dvinskikh2021improved} (статья прошла на конференцию AISTATS 2021), также полученный Дариной, во время проектной смены в Сириусе в августе 2020 года совместно с участником смены -- Д. Тяпкиным 
\begin{center}
  \url{https://sochisirius.ru/news/3901}. 
\end{center}
 В \cite{dvinskikh2021improved} были получены наилучшие известные сейчас оценки временной сложности решения задачи поиска барицентра Васерштейна вероятностных мер. Алгоритм использует не двойственное, а специальное седловое представление исходной задачи. 
 }
\ag{Стоит также отметить, что евклидова регуляризация в задаче \eqref{empiric} нужна, если $\mu < \e/R^2$. При этом до настоящего момента прорабатывалась именно евклидова регуляризация \cite{shalev2009stochastic}. В работе \cite{dvinskikh2020sa} проработан общий случай на примере регуляризации, связанной с 1-нормой (не 2-нормой, как рассмотрено в \eqref{empiric}).}

Без дополнительных предположений о гладкости $f$ в онлайн подходе батч-параллелизация\footnote{Замена в используемых алгоритмах стохастического градиента средним арифметическим стохастических градиентов, рассчитанных параллельно  в одной и той же точке $x$, но с разными (независимыми в совокупности) реализациями.} \cite{gasnikov2020book} возможна, насколько нам известно, только на $\tilde{O}(1)$ ядрах (процессорах) \cite{woodworth2018graph}. При наивном офлайн подходе, в котором для задачи \eqref{empiric} выбирается наилучший способ (метод зеркального спуска \cite{juditsky2012first-order}) решения негладких выпуклых задач (не учитывающий структуру функционала вида суммы) получается оценка $m(\e')$ (здесь и далее в этом абзаце существенно используется то, что $m(\e')$ имеет вид \eqref{SNS} с точностью до логарифмических множителей) для числа необходимых итераций (вычислений $\nabla f(x,\xi^k)$ на каждом узле). После каждой итерации осуществляется коммуникация (всех узлов со всеми или через центр, в зависимости от интерпретации), т.е. всего будет $m(\e')$ актов коммуникации. Ну и число слагаемых в сумме $m(\e)$ также в данном случае отвечает за максимальное число узлов, на котором можно осуществлять полную параллелизацию вычислений. Используя структуру функционала  вида суммы, можно предложить более хитрый способ организации коммуникаций в описанном подходе \cite{lan2017communication,scaman2018optimal,uribe2020dual,dvinskikh2019decentralized}, который приводит к уменьшению числа шагов, на которых необходима коммуникация между узлами (обмен информацией) с $m(\e')$ до $\sqrt{m(\e')}$. Без дополнительных предположений приведенные выше результаты уже не могут быть улучшен ни по одному из показателей \cite{arjevani2015communication,scaman2018optimal}. Из чего можно сделать вывод, что в негладком случае офлайн подход с прямым оракулом сильно уступает онлайн подходу. В более общем случае (децентрализованных подходов) приведенный выше результат отмечен также в таблице~\ref{tab:DetPrimeOr}.

\ag{В заключение этого раздела отметим также недавний результат \cite{gladin2020math}, \cite{gladin2021} (Е. Гладина и К. Зайнуллиной), полученный на проектной смене в Сириусе в августе 2020 \url{http://dmivilensky.ru/opt/}, заключающийся в возможности батч-параллелизации метода Вайды для задачи \eqref{SP} при $\mu = 0$ с размером батча $\tilde{O}\left(\frac{M^2R^2}{{\e}^2}\right)$. При этом число итераций метода будет\footnote{Строго говоря, в этой оценке предполагается компактность множества $Q$, а евклидова асферичность множества $Q$ также должна входить (под логарифмом) в оценку \cite{nemirovsky1983problem,bubeck2014convex}. Аналогичные оговорки следует сделать для офлайн аналога данной формулы, приведенной ниже.} $O\left(n\ln(\Delta f/\e)\right)$. При офлайн подходе возможно распараллеливание метода эллипсоидов для задачи \eqref{empiric} при $\mu = 0$ на $m(\varepsilon)$
процессорах (узлах). При этом число итераций метода будет $O\left(n\ln(\Delta f/\e')\right)$. 
} 
Таким образом, если есть большие возможности в (батч-)параллелизации вычислений, то описанные подходы позволяет существенно ускорить вычисления.

\gav{Кстати говоря, относительно распределенных (децентрализованных) вариантов методов отсечений для задач вида \eqref{empiric} нам практически ничего не известно. Понятно, что с помощью консенсусного алгоритма (см. ниже) и упомянутого анализа чувствительности методов отсечений к неточности субградиента можно, подобно \cite{rogozin2020towards}, получить децентрализованные аналоги методов отсечений. Но не понятно, можно ли в таком классе построить оптимальный алгоритм для задач с небольшим значением $n$? 
}

\subsection{Гладкий случай} \label{section_1.2}
В условиях $L$-Липшицевости градиента $f$ в 2-норме:
\begin{equation}\label{L}
\|\nabla f(y) - \nabla f(x)\|_2 \le L \|y-x\|_2
\end{equation}
и
\begin{equation}\label{sigma}
    \E\left[\nabla f(x, \xi)\right]  \equiv \nabla f(x), \quad \E[\|\nabla f(x,\xi) - \nabla f(x) \|_2^2]\le \sigma^2
\end{equation}
оценки \eqref{SNS} онлайн подхода можно редуцировать до оценок \cite{ghadimi2013stochastic,devolder2013exactness,dvurechensky2016stochastic,gasnikov2018universal}:
\begin{equation}\label{SS}
\hspace{-0.4mm}\min\left\{O\left(\sqrt{\frac{LR^2}{\e}}\right)+O\left(\frac{\sigma^2R^2}{\e^2}\right),~ O\left(\sqrt{\frac{L}{\mu}}\ln\left(\frac{\mu R^2}{\e}\right)\right)+O\left(\frac{\sigma^2}{\mu\e}\right)\right\}.
\end{equation} 
В таких предположениях батч-параллелизация может быть осуществлена на
\begin{equation}\label{batch}
O\left(\frac{{\sigma^2R^2}/{\e^2}}{\sqrt{{LR^2}/{\e}}}\right) ~ \text{ или  } ~ 
O\left(\frac{{\sigma^2}/{(\mu\e)}}{\sqrt{{L}/{\mu}}\ln\left({\mu R^2}/{\e}\right)}\right)
\end{equation}
узлах. Причем, данные оценки в общем случае не могут быть улучшены \cite{woodworth2018graph}. \ag{Более аккуратные (в вероятностном плане) формулировки (включающие проработку возможной неограниченности множества $Q$ -- насколько нам известно, это сделано впервые для класса задач стохастической оптимизации без предположения равномерной ограниченности констант $L$, $\sigma^2$) см. в работе Э. Горбунова и др. \cite{gorbunov2019optimal}. Отметим также работу Д. Двинских, А. Тюрина и др. \cite{dvinskikh2020accelerated}, в которой в максимальной известной нам общности обосновывается тезис о том, что <<оптимальные>> методы для задач стохастической выпуклой оптимизации могут быть построены на базе оптимальных (детерминированных) методов обычной (не стохастической) выпуклой оптимизации за счет правильного выбора размера батча: точность аппроксимации градиента была минимально достаточной, чтобы детерминированный метод <<воспринимал>> пробатченный стох. градиент как настоящий градиент.}

Отметим, что приведенные выше результаты о возможности батч-парал\-лелизации в онлайн режиме естественным образом переносятся и на офлайн режим, поскольку задачу \eqref{empiric} также можно понимать как задачу стохастической оптимизации 
\begin{equation}\label{empiric2}
\min_{x\in Q \subseteq \mathbb{R}^n} \E_k [f(x,\xi^k)] + \frac{\e}{2R^2}\|x - x^0\|_2^2,
\end{equation}
где математическое ожидание берется по случайной величине $k$ равновероятно принимающей значения $1,...,m$.\footnote{На примере такого представления можно пояснить определенную эквивалентность между ранней остановкой процедур обучения типа SGD для задач вида \eqref{empiric}, \eqref{empiric2} (без регуляризации) и ролью регуляризации в офлайн подходе \cite{goodfellow2016deep}. Оказывается \cite{shapiro2005complexity,shalev2009stochastic,shapiro2014lectures,hardt2016train,guigues2017non-asymptotic,lei2020fine}, что если $\mu = 0$ и не добавлять регуляризирующий член, то для $\e$-аппроксимации (по функции) решения задачи \eqref{SP} даже для идеально точного решения задачи \eqref{empiric} или \eqref{empiric2} (без регуляризации) необходимо брать (с точностью до логарифмических множителей) $m \sim nM^2R^2/\e^2$, т.е. в $n$ раз больше, чем с регуляризацией. С другой стороны, если использовать алгоритм SGD для задачи \eqref{empiric2} без регуляризации, то его отличие от SGD для задачи \eqref{SP} будет только в том, что для задачи \eqref{empiric2} без регуляризации выборка независимо осуществляется из заранее зафиксированного набора реализаций объема $m$, в то время как для \eqref{SP} на каждой новой итерации приходит новая независимая реализация $\xi$. Если общее число итераций, которое делает SGD, заметно меньше $m$ (например, для $m \sim nM^2R^2/\e^2$ в теории достаточно $N \sim M^2R^2/\e^2$ итераций; Последующие итерации (в теории) могут только ухудшать качество полученного таким образом решения, т.е. приводить к переобучению; Отметим, что на практике это условие можно существенно ослабить и говорить о сотнях, а иногда даже тысячах <<проходах датасета>>), то различие между тем какая решается задача \eqref{SP} или \eqref{empiric2} без регуляризации не будет особо заметно. Таким образом, ранняя остановка процедуры типа SGD для задачи \eqref{empiric2} без регуляризации приводит к результату аналогичному достаточно точному решению задачи \eqref{empiric} или \eqref{empiric2}. Попытка же точно решить задачу \eqref{empiric}  или \eqref{empiric2} без регуляризации приводит к переобучению. То есть к тому, что полученное таким образом решение будет плохим (по качеству аппроксимации по функции) решением задачи \eqref{SP}. Развитие данного направления частично описано в \cite{sridharan2012learning}. Заметим также, что в последнее время стали достаточно популярны работы (см., например, \cite{mishchenko2020random} и цитированную там литературу) о том, как можно ускорять процедуры типа SGD для задачи \eqref{empiric} за счет случайной перестановки слагаемых в сумме вида \eqref{empiric}, в которой слагаемые уже не обязаны быть реализациями одной и той же случайной функции, и последовательном (по итерациям при проходе всего набора данных) вычислении градиентов соответствующих слагаемых в качестве стохастических градиентов. При таком подходе градиент каждого слагаемого при одном проходе по всему датасету будет использован ровно один раз.

\gav{Стоит отметить, что несмотря на важность затронутых здесь вопросов, на данный момент, насколько нам известно, нет исчерпывающего анализа того, что было описано выше в категориях вероятностей больших отклонений. Наиболее продвинутый результат имеется в работе \cite{feldman2019high}. В невыпуклом случае нам не известно ничего подобного.}} Однако такой способ не приводит к оптимальным оценкам на общее число вычислений стохастических градиентов $f$ и причина не только в том, что задачу \eqref{empiric2} надо решать с более высокой точностью $\e'$, чем исходную задачу. Есть и другая, более важная, причина. Задача вида \eqref{empiric} хотя и может быть представлена в виде \eqref{empiric2}, но все же является достаточно специальным представителем большого класса задач стохастической оптимизации. Именно это обстоятельство позволяет использовать специальный прием <<редукции дисперсии>> (замена стохастического градиента в SGD-типа алгоритмах на некоторый редуцированный стохастический градиент, имеющий тем меньшую дисперсию, чем ближе мы к решению), описанный в следующем разделе. С помощью этого приема удается <<заглянуть>> в структуру задачи и предложить подход, работающий лучше отмеченных выше нижних оценок \cite{woodworth2018graph}.

\ag{Интересно заметить, что если предполагать гладкость (по $x$) $f(x,\xi)$ более высокого порядка чем первого, то подобно отмеченному в предыдущем разделе подходу на базе метода эллипсоидов, задачу \eqref{SP} можно решать тензорными методами \cite{dvinskikh2020accelerated}, что может уменьшить число последовательных итераций.  Результаты, полученные А. Агафоновым, Д. Камзоловым и др. \cite{agafonov2020inexact} при $\mu = 0$ на проектной смене в Сириусе в августе 2020 г. 
\begin{center}
  \url{https://sochisirius.ru/obuchenie/graduates/smena673/3258},  
\end{center}
показывают, что при онлайн подходе есть экономия в размерах батчей для гессиана $\sim \e^{-4/3}$ по сравнению с оценкой размера батча для градиента $\sim \e^{-2}$
для неускоренных тензорных методов второго и третьего порядка \cite{nesterov2018implementable}.\footnote{Используя трюк Ю.Е. Нестерова \cite{nesterov2020inexact,kamzolov2020near}, для методов третьего порядка можно не считать (не оценивать) тензор третьих производных.} При этом число итераций (непараллелизуемых) неускоренного тензорного метода третьего порядка будет $\sim \e^{-1/3}$, что лучше, чем у оптимальных методов первого порядка (градиентных методов) $\sim\e^{-1/2}$  \cite{nesterov2018implementable}. Описанные выше результаты можно распространить на офлайн подход  подобно тому, как это было сделано ранее для метода Вайды.} \gav{Для ускоренных тензорных методов (см., например, \cite{dvinskikh2020meta} и цитированную там литературу) вопрос о получении строгих результатов, насколько нам известно, остается открытым.}

\ag{В ходе работы по проекту данное направление (распределенные тензорные методы) получило развитие в цикле статей \cite{daneshmand2021newton,dvurechensky2021hyperfast,agafonov2021accelerated}.}

\subsection{Редукция дисперсии} \label{section_1.3}
Про метод редукции дисперсии (<<выделения главной части>> по терминологии методов Монте-Карло \cite{ermakov2009method}) написано уже несколько монографий и обзоров, см., например, \cite{lan2020first,linaccelerated,gower2020variance,gasnikov2020book}. \ag{Компактное описание данной конструкции приведено после замечания 1 в приложении в монографии \cite{gasnikov2020book} и в п. E \cite{gower2020variance}.} Далее мы изложим лишь основные результаты, не привязываясь к статистической (вероятностной) специфики возникающих у нас задач минимизации функционалов вида суммы, т.е. далее вместо обозначения $f(x,\xi^k)$ будет использовано обозначение $f_k(x)$. 

Итак, рассмотрим следующую постановку задачи, близкую к \eqref{empiric}, но все же в общем случае отличную от нее (обратим внимание, что здесь $Q = \mathbb{R}^n$, \gav{по-видимому, от этого предположения можно отказаться, однако мы пока не видели, чтобы все это было сделано в такой же общности, в которой проработан случай $Q = \mathbb{R}^n$}):
\begin{equation}\label{empiric3}
\min_{x\in \mathbb{R}^n} f(x) := \frac{1}{m}\sum_{k=1}^m f_k(x),
\end{equation}
где $f_k$ -- выпуклые функции, имеющие $L$-Липшицев градиент в 2-норме, а $f$ -- $\mu$-сильно выпуклая в 2-норме ($\mu \ge 0$).\footnote{Для возможности перенесения обсуждаемых далее результатов на децентрализованные распределенные алгоритмы требуется $\mu$-сильная выпуклость каждой функции $f_k$. Однако от этого обременительного требования легко избавиться, немного видоизменив функции $f_k$ \cite{scaman2017optimal}.} Обратим, внимание, что условие $L$-Липшицев градиента в 2-норме выполняется, например, для квадратичных функций $f_k$. Если переписать задачу \eqref{empiric3} подобно \eqref{empiric2}, то подходы, основанные на различных вариантах SGD, приводящие к оценкам вида \eqref{SNS}, \eqref{SS}, \eqref{batch}, базирующиеся на предположениях \eqref{M}, \eqref{sigma}, уже могут не работать, поскольку отмеченные предположения (равномерно по $x$) уже могут не выполняться (для квадратичных $f_k$). И, действительно, можно даже построить соответствующий пример \cite{assran2020convergence}, когда такая ситуация имеет место. \ag{Это дополнительно указывает на важность анализа, проведенного Э.~Горбуновым и др. в работе \cite{gorbunov2019optimal}, для обоснования оценок \eqref{SS}, \eqref{batch} без предположений равномерного выполнения условия \eqref{sigma}.} Далее, для простоты рассуждений, мы будем исключать такие ситуации, когда один из рассматриваемых (сравниваемых) подходов по каким-то причинам не работает.

Основной результат тут может быть сформулирован следующим образом. Для достижения точности решения  по функции $\e$ для задачи \eqref{empiric3} достаточно 
\begin{equation}\label{VR}
\hspace{-0.4mm}\min\left\{O\left(m+\sqrt{\frac{mLR^2}{\e}}\right),~ O\left(m+\sqrt{\frac{mL}{\mu}}\right)\right\}\ln\left(\frac{\Delta f}{\e}\right)
\end{equation} 
вычислений $\nabla f_k$. Данные результаты в общем случае не могут быть дальше улучшены \cite{woodworth2016tight,lan2020first}. К сожалению, такой метод в общем случае (без дополнительных предположений типа $m > L/\mu$) параллелится лишь на $\tilde{O}(1)$ узлах. В то время как обычный быстрый градиентный метод (см., например, \cite{nesterov2018lectures,gasnikov2020book,lan2020first,linaccelerated}), примененный к задаче \eqref{empiric3} будет параллелиться на $m$ узлах и сходиться за
\begin{equation}\label{N}
N(\e) =  \min\left\{O\left(\sqrt{LR^2/\e}\right),O\left(\sqrt{L/\mu}\ln\left(\mu R^2/\e\right)\right)\right\}
\end{equation}
последовательных итераций \cite{woodworth2018graph}. Причем этот результат будет оптимален в смысле невозможности улучшить число последовательных итераций. Однако, общая трудоемкость такого метода $mN(\e)$ (число вычислений градиентов слагаемых в \eqref{empiric4}) будет в $\sim\sqrt{m}$ раз хуже, чем у подхода с редукцией дисперсии \eqref{VR}.


\begin{itemize}
    \item Прежде всего, заметим, что как онлайн оценки \eqref{SS}, \eqref{batch} так и их офлайн аналог с редукцией дисперсии \eqref{VR} (см. также \eqref{N}),  переходят из сильно выпуклого режима в просто выпуклый при $\mu \sim \e/R^2$, что как раз соответствует выбору коэффициента регуляризации в задаче \eqref{empiric}. Объясняется все это достаточно просто. При таком коэффициенте регуляризации можно заменить исходную задачу \eqref{SP} в случае, если она не сильно выпуклая, на сильно выпуклую, добавив регуляризатор подобно тому, как это сделано в \eqref{empiric}. Тогда $\e/2$-решение (по функции) регуляризованной задачи \eqref{SP} будет $\e$-решением (по функции) исходной задачи \eqref{SP}. При этом выбрать параметр регуляризации больше $\e/(2R^2)$, с сохранением выполнения этого свойства, уже невозможно (см., например, замечание 4.1 \cite{gasnikov2020book}). Собственно, не ограничивая общности, можно было бы и изначально считать, что мы рассматриваем сильно выпуклую задачу с $\mu \ge \e/R^2$. В таком случае регуляризация в \eqref{empiric} будет уже не нужна. Далее в этом разделе мы ограничимся рассмотрением только сильно выпуклого случая.
    \item (\textit{градиентный слайдинг}) В цикле работ Джорджа Лана \cite{lan2020first} (развивающих, в свою очередь, идеи А.С. Немировского и А.Б Юдицкого) была описана конструкция ускоренного градиентного слайдинга. Если говорить совсем грубо, то данная конструкция обосновывает формулу \eqref{VR} при малых $m$. Более точно, эта конструкция говорит, что если сложность решения задачи $\min f(x)$ с точностью $\e$ по функции равна $N_f(\e)$ вычислений $\nabla f$, а сложность решения задачи $\min g(x)$ с точностью $\e$ по функции равна $N_g(\e)$ вычислений $\nabla g$, то сложность решения задачи $\min \{f(x) + g(x)\}$ с точностью $\e$ по функции равна $\tilde{O}(N_f(\e))$ вычислений $\nabla f$ и $\tilde{O}(N_g(\e))$ вычислений $\nabla g$. Дж. Ланом в \cite{lan2020first} было подмечено, что этот прием также можно использовать и для построения распределенных алгоритмов о чем пойдет речь ниже. \ag{В цикле работ авторского коллектива  \cite{beznosikov2019derivative,ivanova2020oracle,gasnikov2020book,dvinskikh2020meta} конструкция слайдинга была заметно расширена. Например, на класс инкрементальных оракулов, рассматриваемых в данном разделе, покомпонентных и безградиентных оракулов; в том числе и смешанных -- по $f$ имеем $\nabla f$, а по $g$ имеем доступ только к значениям $g$ \cite{beznosikov2019derivative}. Такие обобщения во многом также были мотивированы разработкой оптимальных распределенных децентрализованных алгоритмов, см. \cite{dvinskikh2019decentralized,beznosikov2019derivative} и следующий раздел.}
    \item Если в \eqref{VR} подставить оценку $m$ из \eqref{SNS} и считать, что $L < M^2/\e$ (если это условие не выполняется, то нет смысла считать задачу гладкой -- оценка \eqref{SS} будут хуже оценки \eqref{SNS}), то оценка \eqref{VR} с точностью до логарифмических множителей будет просто совпадать с $m$. Поскольку для задачи \eqref{empiric} в оценке \eqref{VR} минимум достигается на втором аргументе, то точность $\e'$ решения задачи \eqref{empiric} будет входить в оценку сложности \eqref{VR} только под логарифмом. Таким образом, 
    сложность офлайн подхода на базе подхода с редукцией дисперсии  совпадает (также с точностью до логарифмических множителей) с точностью онлайн подхода, что <<восстанавливает справедливость>>, которая, как могло показаться, была потеряна в предыдущем разделе при обсуждении тех сложностей, которые возникают при решении задачи \eqref{empiric2} без редукции дисперсии.   
\end{itemize}

Наряду с задачей \eqref{empiric3} можно рассмотреть задачу
\begin{equation}\label{empiric4}
\min_{x\in \mathbb{R}^n} f(x) := \frac{1}{m}\sum_{k=1}^m f_k(x) = \frac{1}{m}\sum_{k=1}^m \EE[f_k(x,\xi^k)],
\end{equation}
считая доступным только стохастические градиенты слагаемых, удовлетворяющие для любого $k = 1,...,m$ условию \eqref{sigma}. Для такой задачи оценка \eqref{VR} изменится следующим образом \cite{kulunchakov2019estimate1,kulunchakov2019estimate2,kulunchakov2019generic}:
\begin{equation}\label{SVR}
    m+\sqrt{\frac{mL}{\mu}} \to m+\sqrt{\frac{mL}{\mu}} + \frac{\sigma^2}{\mu\e}
\end{equation}
и аналогично в выпуклом случае: в \eqref{SVR} нужно подставить $\mu \sim \e/R^2$, см. выше. При этом за счет батчинга можно полностью редуцировать последнее слагаемое подобно \eqref{batch}. 

Если вместо редукции дисперсии использовать обычный ускоренный градиентный метод с помощью которого для задачи \eqref{SP} (или что то же самое \eqref{empiric4} с $m = 1$) были получены оценки \eqref{SS}, \eqref{batch}, то оценка \eqref{SVR} перейдет в оценку 
\begin{equation}\label{SVRFGM}
    m+\sqrt{\frac{mL}{\mu}} + \frac{\sigma^2}{\mu\e} \to m\cdot\left(\sqrt{\frac{L}{\mu}} +  \frac{\sigma^2}{m\mu\e}\right),
\end{equation}
которая очевидным образом уже параллелится на $m$ узлах (при необходимости последнее слагаемое за счет батчинга может быть дополнительно редуцировано). Без дополнительных предположений оценка \eqref{SVRFGM} не может быть улучшена \cite{woodworth2018graph,woodworth2020minibatch}. Точнее говоря, не может быть улучшена ни оценка на число последовательных итераций / коммуникаций (получается за счет батчинга) $\sim\sqrt{L/\mu}$, ни оценка числа вычислений $\nabla f_k(x,\xi^k)$ на каждом из $m$ узлов (аналогично \eqref{SS}).

Если же вместо $f_k(x) = \E[f_k(x,\xi^k)]$ в \eqref{empiric3} стоит $f_k(x):=\frac{1}{r}\sum_{j=1}^r f_k^j(x)$, где все $f_k^j$ удовлетворят условию \eqref{L},  то оценка \eqref{SVRFGM} может быть уточнена следующим образом \cite{li2020optimal,gorbunov2020local,rajawat2020primal}:
\begin{equation}\label{VRVR}
    m+\sqrt{\frac{mL}{\mu}} + \frac{\sigma^2}{\mu\e} \to m\cdot\left(r+\sqrt{\frac{rL}{\mu}}\right) 
\end{equation}
вычислений $\nabla f_k^j$.
Причем здесь, как и в случае с ускоренным градиентным методом, работающим по оценке \eqref{SVRFGM}, возможна параллелизация на $m$ узлах, приводящая к одной и той же оценке числа числа шагов с обменом информацией между узлами, т.е. оценке числа коммуникаций (см. также \eqref{N} и \eqref{SS}, \eqref{batch}): $\sim\sqrt{L/\mu}$. Эта оценка вместе с оценкой \eqref{VRVR} не могут быть улучшены \cite{hendrikx2020optimal}. Неулучшаемость оценки $\sim\sqrt{L/\mu}$ означает, что без дополнительных предположений не существует способа решить задачу на $m$ узлах с числом шагов коммуникаций, на которых осуществляется обмен информацией между узлами, меньше, чем  $\sim\sqrt{L/\mu}$. Неулучшаемость оценки \eqref{VRVR} понимается в том смысле, что для любых параллельных/распределенных алгоритмов число вычислений $\nabla f_k^j$ на каждом из $m$ узлов не может быть в общем случае сделано меньше, чем $\sim r+\sqrt{rL/\mu}$. Аналогичный результат мы уже упоминали ранее при обсуждении формулы \eqref{VR} \cite{woodworth2016tight,lan2020first}.

Результаты \eqref{SVRFGM}  и \eqref{VRVR} можно сравнить на задаче \eqref{empiric} c $m$ определяемым \eqref{SNS}. Редуцируя исходную сумму из $m$ слагаемых в сумму $m/r$ слагаемых, каждое из которых представляет собой, в свою очередь, сумму $r$ слагаемых (при подходе, приводящем к формулам \eqref{SVRFGM}, эта сумма $r$ слагаемых представляется в виде математического ожидания по равномерной мере подобно представлению \eqref{empiric2}), получим, что при одинаковом числе коммуникаций подходы будут приводить к различным оценкам оракульной сложности на узлах. А именно, для подхода, приводящего к оценке \eqref{VRVR}, число вызовов оракула в каждом узле будет  $\tilde{O}\left(\left(r+\sqrt{rL/\mu}\right)\ln\left(\Delta f/\e'\right)\right)$, а для подхода, приводящего к оценке \eqref{SVRFGM}, оценка будет\footnote{Считаем здесь и далее (таблица~\ref{T:stoch_prima_oracle}) $\sigma^2 \sim M^2$.} $O\left(\sqrt{L/\mu} + r\e/\e'\right)$. Учитывая, что точность решения задачи \eqref{empiric} $\e' = O\left(\mu\e^2/M^2\right)$, получаем, что подход базирующийся на редукции дисперсии и приводящий к оценке $\sim r+\sqrt{rL/\mu}$ является предпочтительнее.

Последний результат является, пожалуй, наиболее важным наблюдением данного раздела в контексте приложения различных вариаций на тему редукции дисперсии к практическому решению больших задач, приходящих из анализа данных, рассмотренных в предыдущих разделах. Особенно ценно, что этот подход допускает распределенную децентрализованную версию. Изложению децентрализованных версий рассматриваемых до сих пор подходов и посвящен следующий раздел.

\gav{Главным же недостатком описанного в данном разделе (и предыдущих разделах) офлайн подхода(-ов), является предположение, что функции $f_k$ разные. На самом деле во многих реальных приложениях из анализа данных $f_k$ статистически близки, поскольку являются, в свою очередь, суммами одинаково распределенных случайных функций, хранящихся на соответствующих узлах. Попытки использовать эту статистическую близость представляют собой, на наш взгляд, наиболее перспективные исследования в данной области \cite{arjevani2015communication}. Отметим в качестве примера недавнюю работу \cite{hendrikx2020statistically}, в которой рассматривается централизованная распределенная архитектура коммуникационной сети. Насколько нам известно, какой-либо законченной общей теории здесь пока еще не построено. Тем не менее, отметим работу \cite{sun2020convergence}, в которой демонстрируется как можно использовать статистическую близость в децентрализованной распределенной архитектуре (в случае неускоренных методов и без редукции дисперсии). А именно, происходит редукция $L$: $L \sim \mu+ \frac{\text{const}}{\sqrt{r}}$ в оценках числа шагов коммуникаций, что особенно актуально в распределенном контексте (см. следующие разделы). С построением ускоренных методов возникают сложности обсуждаемые в работе \cite{hendrikx2020statistically}. Тем не менее, по-видимому, сочетание \cite{hendrikx2020statistically,stonyakin2020inexact,rogozin2020towards} может позволить частично перенести результаты \cite{hendrikx2020statistically}, содержащие специальным образом ускоренные методы (в централизованной архитектуре), на общие децентрализованные схемы.} \ag{В ходе работы по проекту в этом направлении были достигнуты определенные успехи \cite{daneshmand2021newton,dvurechensky2021hyperfast,agafonov2021accelerated}. В частности, были построены ускоренные методы на базе методов второго порядка \cite{agafonov2021accelerated}.}

Отметим, что в случае $f_k(x):=\frac{1}{r}\sum_{j=1}^r f_k^j(x)$ для задачи \eqref{empiric3} можно использовать как описанный выше подход с редукцией дисперсии, приводящей к оценке \eqref{VRVR}, так и (ускоренный) градиентный метод с $L \sim \mu+ \frac{\text{const}}{\sqrt{r}}$. Если не учитывать возможности дополнительной параллелизации при вычислении градиента на каждом узле при втором подходе, то выигрыш второго подхода по числу коммуникаций (итераций) будет происходить на фоне проигрыша в числе вычислений $\nabla f_k^j$ на каждом узле. \gav{Возможно ли как-то сочетать редукцию дисперсии в смысле \eqref{VRVR} и  статистическую похожесть функций $f_k$? -- насколько нам известно, открытый вопрос. Также открытым остаются вопросы о распространении обсуждаемых здесь результатов возможности учета статистической близости слагаемых на случай негладкой функции $f$ (см. предыдущий раздел). По-видимому, здесь может помочь регуляризация прямой задачи и переход к двойственной, которая в этом случае уже будет гладкой \cite{uribe2020dual,dvinskikh2019decentralized}.}

\section{Распределенные алгоритмы и \ag{основные полученные результаты}} \label{section_2}
Прежде чем переходить к изложению общих результатов по децентрализованным распределенным алгоритмам напомним (резюмируем) то, что было написано выше для параллельных $\chi = 1$ и распределенных на $m$ узлах централизованных ($\sqrt{\chi} = d$, где $d$ -- диаметр коммуникационной сети \cite{scaman2017optimal}) алгоритмах. Рассматривались задачи вида 
\begin{equation}\label{empiric44}
\min_{x\in \mathbb{R}^n} f(x) := \frac{1}{m}\sum_{k=1}^m f_k(x) = \frac{1}{m}\sum_{k=1}^m \E[f_k(x,\xi^k)]
\end{equation}
и
\begin{equation}\label{empiric5}
\min_{x\in \mathbb{R}^n} f(x) := \frac{1}{m}\sum_{k=1}^m f_k(x) = \frac{1}{m}\sum_{k=1}^m \sum_{j=1}^r f_k^j(x).
\end{equation}

Для задачи \eqref{empiric44} имеющиеся результаты собраны в таблицах \ref{tab:DetPrimeOr} и \ref{T:stoch_prima_oracle}, позаимствованных (кроме результатов, выделенных \agv{красным}) из \cite{dvinskikh2019decentralized} (почти все результаты из таблицы~\ref{tab:DetPrimeOr} были получены в \cite{arjevani2015communication,scaman2017optimal,lan2017communication}).

Для задачи \eqref{empiric5}  число коммуникационных раундов будет $\tilde{O}(\sqrt{\chi L/\mu})$, а число вызовов оракула $\nabla f_k^j$ на каждом узле будет $\tilde{O}\left(r+\sqrt{rL/\mu}\right)$ \cite{li2020optimal}.

	\begin{table}[H]
\caption {Оптимальные оценки для детерминированного оракула $\nabla f_k$}
\label{tab:DetPrimeOr}
\begin{center}
\begin{tabular}{|c| c| c| c| c|}
 \hline
 & \makecell{ $f_k$ is \\ $\mu$-str. convex\\ and $L$-smooth} 
 &   \makecell{$f_k$ is $L$-smooth} 
 & \makecell{ $f_k$ is \\$\mu$-str. convex} &  \\
 \hline
 \makecell{\#communic. \\ rounds} 
 & \makecell{$\widetilde O\left(\sqrt{\frac{L}{\mu}\chi} \right)$}   
 & \makecell{ $ \widetilde O\left({\sqrt{\frac{LR^2}{\e} \chi}}\right)$ }
 &  \makecell{  $O\left(\sqrt{\frac{M^2}{\mu\e}\chi} \right)$} 
 &  \makecell{  $O\left( \sqrt{\frac{M^2R^2}{\e^2}\chi} \right) $}   \\
 \hline
\makecell{\#oracle calls\\of $\nabla f_k$\\ per node $k$} 
& $ \widetilde O\left(\sqrt{\frac{L}{\mu}} \right)$ 
&  $O\left(\sqrt{\frac{LR^2}{\e}}  \right) $ 
&   $O\left(\frac{M^2}{\mu\e} \right) $ 
&   $O\left( \frac{M^2R^2}{\e^2}\right)$\\
 \hline
Algorithm & \makecell{{\tt PSTM}, \\ $Q=\R^n$} & \makecell{{\tt PSTM},\\ $Q=\R^n$} &  {\tt R}-Sliding &  Sliding\\
 \hline
\end{tabular}
\end{center}
\end{table}

	\begin{table}[H]
\caption {Оптимальные оценки для стохастического оракула $\nabla f_k(x,\xi^k)$} 
\label{T:stoch_prima_oracle}
\begin{center}
\begin{tabular}{ |c| c| c| c| c|}
 \hline
&{\makecell{ $f_k$ is \\
$\mu${-str. convex}\\ and $L$-smooth} } &  {\makecell{ $f_k$ is \\ $L$-smooth} }&  {\makecell{ $f_k$ is \\
$\mu${-str. convex}} } &  \\
 \hline
\makecell{\#communic. \\
rounds}& \makecell{ $\widetilde O\left(\sqrt{\frac{L}{\mu}\chi} \right)$}   & \makecell{ $ \widetilde O\left({\sqrt{\frac{LR^2}{\e} \chi}}\right)$ }&  \makecell{  $O\left(\sqrt{\frac{M^2}{\mu\e}\chi} \right)$} &  \makecell{  $O\left( \sqrt{\frac{M^2R^2}{\e^2}\chi} \right) $}  \\
 \hline
   \makecell{\#oracle calls \\ of $\nabla f_k(x, \xi^k)$\\
 per node $k$} & {\makecell{ \agv{$\widetilde O\left(\max\left\{\frac{\sigma^2}{m\mu\e},
\right.\right.$}\\
\agv{$\left.\left.  \sqrt{\frac{L}{\mu}}\right\}\right)$} }}  & {\makecell{\agv{$ O\left(\max\left\{\frac{\sigma^2R^2}{m\e^2}, \right.\right.$}\\ 
\agv{$\left.\left.\sqrt{\frac{LR^2}{\e}}  \right\}\right)$}}} & { {\makecell{  $O\left(\frac{M^2}{\mu\e} \right)$}}} &{\makecell{  $O\left(\frac{M^2R^2}{\e^2}\right)$}} \\
 \hline
Algorithm & \makecell{{\tt \cite{dvinskikh2020accelerated,rogozin2020towards}},\\ $Q=\R^n$}& \makecell{{\tt \cite{dvinskikh2020accelerated,rogozin2020towards}},\\ $Q=\R^n$} & \makecell{ Stochastic \\{\tt R}-Sliding} & \makecell{ Stochastic \\ Sliding} \\
 \hline
\end{tabular}
\end{center}
\end{table}

Все приведенные здесь результаты, как уже отмечалось ранее, оптимальны, т.е. не могут быть улучшены без дополнительных предположений (<<дополнительного заглядывания в структуру задачи>>).

\gav{Отметим, что современные исследования по распределенным архитектурам типа федеративного обучения \cite{woodworth2020minibatch} (чередования нескольких последовательных итераций на узлах с шагами коммуникаций всех со всеми, ну или через центр в зависимости от интерпретации) показывают, что для разных гладких функций $f_k(x)=\E[f_k(x,\xi^k)]$ в общем случае (без наличия представлений $f_k$ в виде суммы \cite{gorbunov2020local} и других уточняющих предположений) возможность использовать последовательные итерации на узлах без коммуникации ничего не дает (достаточно батч-параллелизации). По мере того, что функции  $f_k$ становятся близкими друг другу (в пределе одинаковыми \cite{godichon2020rates,woodworth2020local}) появляется возможность дополнительного (к батч-параллелизации) ускорения за счет последовательных итераций на узлах, что может в итоге существенно экономить число коммуникационных шагов. Вплоть до необходимости осуществления всего одного коммуникационного шага в самом конце \cite{godichon2020rates}.
Законченной теории, насколько нам известно, здесь пока еще нет. Это направление представляется, на наш взгляд, одним из самых интересных в современной распределенной оптимизации. Недавняя работа \cite{woodworth2021min} существенным образом проливает свет на то, как будет выглядеть эта теория.}

\gav{Приведем в заключение этого раздела результат из работы \cite{woodworth2020local} для задачи \eqref{empiric44} с $\mu = 0$, в которой все $f_k$ равны между собой (и являются квадратичными функциями). Обозначим через $K$ -- число коммуникаций (на каждой коммуникации происходит обмен сообщениями всех узлов со всеми или с центром, в зависимости от интерпретации), $T$ -- число последовательных итераций (между двумя коммуникациями) на каждом узле (на каждой итерации можно один раз посчитать $\nabla f_k(x,\xi^k)$), $m$ -- число узлов. Тогда оптимальный алгоритм
должен выдавать такой $\tilde{x}$, что\footnote{Эта оценка соответствует левой оценке в \eqref{SS}, т.е. может интерпретироваться, как результат работы после $KT$ итераций ускоренного градиентного метода с батчем размера $m$.} 
\begin{equation}\label{FL}
    \EE [f(\tilde{x})] - f(x_*) \simeq \frac{LR^2}{(KT)^2} + \frac{\sigma R}{\sqrt{mKT}}. 
\end{equation}
Отсюда будет следовать (имеется аналогия с формулой \eqref{batch}) возможность параллелизации вычислений на $m \sim N^{3/4}$ узлах, где $N = mKT$ -- общее число вычислений $\nabla f_k(x,\xi^k)$.
В случае не квадратичных функций оценка \eqref{FL} испортится \cite{woodworth2021min}. Отметим, что имеются алгоритмы, которые имеют в \eqref{FL} вместо $\frac{LR^2}{(KT)^2}$ неускоренную сходимость $\frac{LR^2}{KT}$, что приводит к более скромной оценке $m \sim N^{1/2}$ \cite{godichon2020rates} (в этой работе рассматривается сильно выпуклый случай с $K = 1$). Отметим, что в случае разных функций $f_k$ в той же самой архитектуре оптимальный алгоритм работал бы по формуле аналогичной \eqref{FL} с заменой    $\frac{LR^2}{(KT)^2}$ на $\frac{LR^2}{K^2}$, что альтернативным образом (к последовательным итерациям) можно обеспечить за счет батчинга на каждом узле $\sigma^2 \to \sigma^2/K$.
}

\subsection{Децентрализованная распределенная оптимизация} \label{section_2.1} 
В децентрализованной распределенной оптимизации вводится понятие (связанного) коммуникационного графа на $m$ вершинах. В $k$-м узле хранится $f_k$, точнее говоря, есть возможность обращаться к оракулу, выдающему определенную информацию об $f_k$. В отличие от централизованной архитектуры, в которой выделяется центральный узел, собирающий с остальных узлов-исполнителей информацию, осуществляющий вычисления и распространяющий новую информацию обратно, в децентрализованной оптимизации (всем узлам) разрешены за один раунд коммуникации только со своими непосредственными соседями. Централизованный сбор информации не осуществляется. Если сопоставить коммуникационному графу матрицу Лапласа (являющуюся неотрицательно определенной), то отношение максимального собственного значения этой матрицы к минимальному неотрицательному $\chi$ будет определять (ускоренный консенсусный алгоритм -- см., например, текст после упражнения 4.8 \cite{gasnikov2020book} и \cite{uribe2020dual}) время $\tilde{O}(\sqrt{\chi})$, необходимое для достижения консенсуса на таком графе, т.е. время (а точнее, число коммуникационных раундов), необходимое узлам, чтобы узнать среднее арифметическое чисел, изначально записанных в этих узлах. Эта величина $\tilde{O}(\sqrt{\chi})$ является аналогом диаметра графа $d$. Точнее, она оценивается снизу диаметром графа. Но типично она и равна диаметру графа с точностью до логарифмических множителей. Впрочем, есть ситуации, например, звездная централизованная архитектура коммуникационной сети, когда разница может быть в $m$ раз, и даже больше \cite{nedic2019graph}. \ag{Собственно, можно было бы ожидать, что если в приведенных в предыдущем разделе оценках под $\chi$ понимать то, что было определено в этом разделе, то все результаты (возможно, с небольшими оговорками) удастся сохранить. Так оно и есть на самом деле. Примечательно, что хотя данное направление имеет достаточно длительную историю, см. работы Бертесекаса--Цициклиса, А. Недич и др., лишь в последние 5 лет оно приобрело описываемый здесь вид.\footnote{Особенно отметим работы \cite{arjevani2015communication,shi2015extra,nedic2017achieving,lan2017communication,scaman2017optimal} с хорошим запасом оригинальных идей. С этих статей началась новая волна теоретического интереса к выпуклой децентрализованной распределенной оптимизации. Эта волна не прошла до сих пор.} Обратим внимание, что таблица \ref{T:stoch_prima_oracle} была сформирована во многом по разработкам, полученным А. Гасниковым,  Э. Горбуновым, Д. Двинских, А. Рогозиным и др. в ходе работы по проекту, см. \cite{dvinskikh2019decentralized,gorbunov2019optimal,rogozin2020towards}. Красным цветом в таблице \ref{T:stoch_prima_oracle} выделены результаты, которые получены совсем недавно \cite{rogozin2021accelerated}.  В основе подхода, позволяющего это сделать, лежит консенсусный вариант ускоренного метода из работы \cite{rogozin2020towards}, в который следует внести специальный батчинг \cite{gorbunov2019optimal,dvinskikh2020accelerated} и заменить неускоренную процедуру консенсуса на ускоренную. Отметим также, что в работах \cite{dvinskikh2019decentralized,gorbunov2019optimal,uribe2020dual} (подготовка статьи \cite{uribe2020dual}, опубликованной в 2020 г. в журнале Optimization Methods and Software,  завершалась в ходе работы над проектом) аналогичные таблицы были построены и для двойственного оракула, выдающего (стохастический) (суб-)градиент двойственной функции $f_k^*$. Эти результаты, как уже отмечалось, использовались для построения эффективных алгоритмов поиска барицентра Васерштейна вероятностных мер \cite{dvurechenskii2018decentralize,dvinskikh2019dual,kroshnin2019complexity}.}

\ag{Отметим, что до сих пор мы в основном говорили об оракуле первого порядка, выдающем (стохастический) (суб-)градиент $f_k$. На самом деле, по целому ряду причин градиентный оракул может быть недоступен \cite{spall2003introduction}. Тогда можно использовать безградиентный оракул, выдающий значение (реализацию) функции $f_k$. Естественная идея, активно развиваемая в работах К. Шайнберг и др., что можно просто восстанавливать (приближенно) градиент по конечным разностям. Конечно, в негладком случае при таком подходе могут возникать некоторые сложности, требующие определенных трюков \cite{shamir2017optimal,bayandina2018gradient}, но все же главный вопрос, который тут возникает: возможно ли получить что-то лучше, чем при описанном подходе, который, очевидным образом, редуцирует случай безградиентного оракула к рассмотренному градиентному? Вопрос возник потому, что не в распределенной оптимизации это возможно \cite{duchi2015optimal,gasnikov2016gradient-free,gasnikov2017stochastic,shamir2017optimal,bayandina2018gradient,gorbunov2018accelerated,dvurechensky2020accelerated}. В негладком случае это оказалось возможным и для децентрализованной (стохастической) распределенной оптимизации. Оптимальные алгоритмы были разработаны А. Безносиковым, Э. Горбуновым и А. Гасниковым в статье на IFAC 2020 \cite{beznosikov2020derivative} на основе специальной конструкции: упомянутого ранее слайдинга Лана, использованного в варианте для двух градиентных оракулов, один из которых выдает градиент, а второй стохастический субградиент, который формируется на базе специальной \cite{shamir2017optimal} безградиентной аппроксимации (стохастического) субградиента.} \ag{Результаты в гладком случае, которые развивали бы работы \cite{gorbunov2018accelerated,dvurechensky2020accelerated} нам не известны.}

\subsection{Возможное развитие и обобщения} \label{section_2.2} 
Ранее в ходе изложения мы периодически обсуждали открытые вопросы. Настоящий раздел преимущественно весь состоит из таких обсуждений.

\begin{itemize}
    \item \ag{В работе \cite{parsegov2019accelerated}, принятой на 3rd IFAC Workshop on Cyber-Physical and Human Systems, было продемонстрировано, что если в постановке задачи \eqref{empiric4} аргументы у функций $f_k$ не обязательно одинаковые, а состоят из блоков, которые уникальны для данного слагаемого и блоков, которые встречаются и у других слагаемых, то такую задачу также можно решать распределенными алгоритмами. Было бы интересно перенести описанные выше результаты на постановку задачи из статьи \cite{parsegov2019accelerated}.}
    \item \ag{В начале отчета мы уже упоминали о статье Д. Двинских и Д. Тяпкина \cite{dvinskikh2020improved}, в которой была предложена седловая переформулировка задачи поиска барицентра Васерштейна вероятностных мер с носителем на $n$ точках. Предложенный алгоритм в теории в $\sqrt{n}$ раз работает лучше известных ранее. Более того, есть основания полагать (по аналогии с \cite{blanchet2018towards}), что полученная в \cite{dvinskikh2021improved} оценка трудоемкости подхода уже не может быть улучшена (другими алгоритмами) более чем на логарифмический множитель.} Однако алгоритм Двинских--Тяпкина был описан в не распределенном варианте, хотя задача сформулирована как седловая с функционалом вида суммы, в которой слагаемые имеют как уникальные блоки переменных, так и общий блок переменных (один у всех). Подобно уже отмеченному выше обобщению работы \cite{parsegov2019accelerated} представляется интересным перенести описанные в данном отчете результаты на указанные (выпукло-вогнутые) седловые постановки задач. Ожидается, что для получения нужного обобщения потребуется распространить градиентный слайдинг на седловые задачи, что можно сделать на базе конструкции из п. 3.4 статьи \cite{dvinskikh2020meta}. \ag{Недавно это было сделано в цикле статей \cite{beznosikov2020local,rogozin2021decentralized}.}
    \item В работе большого авторского коллектива из EPFL А. Колосковой, С.~Стича, М. Ягги и др. \cite{koloskova2020unified}, активно работающего в направлении развития различных распределенных алгоримтов, был описан общий подход, объединяющий всевозможные распределенные схемы решения задач выпуклой оптимизации. Это подход включает централизованную и децентрализованную оптимизацию, схему федеративного обучения, распределенную оптимизацию на меняющихся со временем коммуникационных графах, госсип и многое другое. Полученные результаты были сформулированы подобно тому, как были сформированы таблицы \ref{tab:DetPrimeOr}, \ref{T:stoch_prima_oracle}. Однако все это было сделано в неускоренном случае -- для конкретного метода <<распределенного SGD>>. Этот метод не дает оптимальные оценки. Естественно, было в общности статьи \cite{koloskova2020unified} получить описанные выше результаты. \ag{Процесс уже начался. Так в работе \cite{rogozin2020towards} результаты А.~Колосковой и др. ускоряются в случае детерминированного оракула и детерминировано меняющейся со временем коммуникационной матрицы. А в работе \cite{rogozin2021accelerated} это сделано для задач стохастической оптимизации. Отметим, что по консенсусной части ускорения нет, потому что допускается, что граф меняется со временем.} В статье \cite{koloskova2020unified} есть также ограничение, что коммуникационная матрица дважды стохастическая. Это предположение было введено в статье, чтобы единообразно представить большое число результатов. В действительности, заметная часть результатов, по-видимому, может быть получена и без этого предположения.
    \item Важно отметить, что выше, чтобы упростить изложение, мы работали с одной константой Липшица градиента $L = L_f$ у $f$ и $L = L_{f_k}$ у $f_k$. На самом деле, у слагаемых $f_k$ могут быть разные константы $L_{f_k}$, и в разные приводимые выше оценки,  могут входить разные константы $L_f$, $\bar{L} = m^{-1}\sum_{k=1}^m L_{f_k}$, $L_{\max} = \max_k L_{f_k}$, которые могут отличаться в $n$ раз \cite{hendrikx2020optimal,gasnikov2020book} за счет различий у $L_{f_k}$ между собой и различий с $L = L_f$. \ag{В работах \cite{ye2020multi,rogozin2020towards,li2020optimal,rogozin2021decentralized} (статьи \cite{rogozin2020towards,rogozin2021decentralized} были подготовлены в рамках работы по проекту) для детерминированных оракулов исследовался данный вопрос.} Однако на данный момент тут имеются лишь частные результаты. Общего понимания о том где и какие константы гладкости (средние, худшие) стоит писать пока нет. Не известен в общем случае и ответ на вопрос: можно ли как-то уходить от худших константа за счет дополнительных рандомизаций и т.п.? 
    \item Интересные результаты могут быть получены, при обобщении распределенных алгоритмов на случай, когда на передаваемые векторы  (например, стохастические градиенты) перед отправкой действуют оператором квантизации $Q(z)$ (сжатия): $$\E\left[\|Q(z) - z\|_2^2\right]\le (1 - q) \|z\|_2^2.$$
    Примеры таких операторов \textit{TopK} (обнуляет все компоненты вектора, кроме $K$ наибольших по модулю) или \textit{RandK} (обнуляет все компоненты, кроме случайно выбранных $K$ компонент и осуществляет масштабирование с множителем $n/K$). 
    Для этих операторов $q = K/n$. \ag{Естественно ожидать, что многие приведенные выше результаты переносятся на протоколы со сжатием передаваемой информации так, что число коммуникаций возрастает приблизительно в $1/q$ раз \cite{beznosikov2020biased,qian2020error,albasyoni2020optimal,gorbunov2020linearly}.} 
    Интересно понять, насколько этот результат является общим и что происходит с оракульной сложностью? \ag{В работе \cite{krawtschenko2020distributed} был предложен рандомизированный способ сжатия данных при передаче, перекликающийся с техникой рандомизации Григориадиса--Хачияна \cite{grigoriadis1995sublinear}. А именно, вместо вектора распределения вероятностей передается одна из случайно разыгранных (согласно этому вектору) вершин симплекса. Такой способ оказался достаточно эффективным при решении двойственной задачи для задачи поиска барицентра Васерштейна вероятностных мер. Стохастический градиент двойственной функции к энтропийно-сглаженному расстоянию Васерштейна как раз оказался вектором из единичного симплекса (распределением вероятностей). Отметим также, что такая квантизация выполняет и другую функцию -- обеспечивает большую конфиденциальность передаваемых данных. }
\end{itemize}

Авторский коллектив выражает благодарность А.С. Ненашеву за организацию августовских проектных смен в Сириусе (г. Сочи), на которых была получена (доработана) заметная часть вошедших в отчет результатов, а также за постоянную помощь по широкому спектру вопросов во время этих смен.

\bibliographystyle{abbrv}
\bibliography{PD_references}

\end{document}